\documentclass[12pt,a4paper]{article}

\usepackage[english]{babel}
\usepackage[utf8x]{inputenc}
\usepackage{amsmath}
\usepackage{amssymb}
\usepackage{graphicx}
\usepackage{amsthm} 

\usepackage[colorlinks=true, allcolors=blue]{hyperref}
\usepackage{url}

\newtheorem{theorem}{Theorem}[section]

\newtheorem{lemma}{Lemma}[section]
\newtheorem{assumption}{Assumption}[section]

\title{Oscillation time and damping coefficients  in a nonlinear pendulum}
\author{Jaime Arango}
\begin{document}
\maketitle

\begin{abstract} 
We establish a relationship between the normalized damping coefficients and the time that takes a nonlinear pendulum to complete  one oscillation   starting from an initial position with vanishing velocity. We establish some conditions on the nonlinear restitution force  so that    this oscillation time    does not depend monotonically on the viscosity damping  coefficient. 

\smallskip
\noindent \textbf{ASC2020:} 34C15, 34C25

\smallskip
\noindent \textbf{Keywords.} oscillation time, damping, damped oscillations

\smallskip
\begin{quotation}
This paper is dedicated to the memory of Prof. Alan Lazer (1938-2020), University of Miami. It was my  pleasure to discuss with him some of the results presented here
\end{quotation}
\end{abstract}

\section{Introduction}
The pendulum is perhaps the oldest and fruitful paradigm for the study of an oscillating system. The apparent regularity of an oscillating mass going to and fro  through the equilibrium position has fascinated the scientists well before Galileo.  There are plenty of  mathematical models accounting for almost any observed behavior of the  pendulum's  oscillation. From the sheer amount of the literature on the subject,  one would expect that there is no reasonable question regarding a pendulum that has no been already answered. And that might be true. Yet, for  whatever reason,  it is not impossible to take on a question whose answer does not seem to follow immediately from the classical sources. 

In a typical experimental setup with no noticeable damping, the oscillations of a pendulum are periodic. Now,  if the damping cannot be neglected, we still observe oscillations, even though they are non periodic. However, we can measure the time spent by a complete oscillation, and this time is a natural generalization of the period.  But, how does depend this oscillation time on the characteristic  of the medium, say on the viscosity of the surrounding atmosphere? It seems that there is no much information on how the damping affects the oscillation time. There are plenty of new publications regarding damping and  oscillations,  ranging from analytical solutions (\cite{Johannessen_2014}, \cite{ghosechoudhury},\cite{kharkongor}), to very clever experimental setups (see for example  \cite{marmolejo}). The nature of the damping has been also extensively considered (\cite{Zonetti_1999}, \cite{cveticanin}), but the dependence of the oscillation time on the damping or on the non-linearity seems to be less investigated.

For the sake of simplicity we analyze the oscillation time  in the frame of a model that appear in almost any text book of ordinary differential equations  (see for example \cite{arnold1}):   
\begin{equation}
\label{eq:pendulum}
\ddot{x} +2\alpha\, \dot{x} +x\left( 1+ f\left(x\right)\right)= 0,
\end{equation}
where $x = x(t)$ measures  the pendulum's deviation   with respect to a vertical axis of equilibrium and  $\alpha \ge 0$ denote the viscous damping coefficient. The term $x\, f(x)$ models the nonlinear  part of the restoring force. We've rescaled the time so that  the period of the linear undamped oscillation is exactly $2\pi$. 

The math of the solutions $x = x(t)$ is classical. If $f$ is smooth and $x_0$ and $v_0$ are given real values, then there exists a  unique solution  satisfying the given conditions $x(0)= x_0$ and $ \dot{x}(0) = v_0.$ Moreover, if $f(0)= 0$,  then $x = 0$ is a  stable equilibrium solution of \eqref{eq:pendulum}. As a consequence, $x(t)$ is defined for all $t\ge 0$ provided $|x_0|\ll 1$ and $|v_0|\ll 1$. Notice that the points of  vanishing derivative of a solution $x= x(t)$ to \eqref{eq:pendulum} are isolated and those points correspond,  either to local maxima or to local minima. Denote by $\tau(x_0,\alpha)$ the amount of time spent (by the mass) completing  one oscillation starting  from $x_0$ with vanishing velocity  ($v_0 = 0$). To be precise, if $x = x(t)$ starts from $x_0$ with vanishing velocity, then $x$ reaches a local maximum at $t = 0$, and  the oscillation is completed when $x$ reaches the next local maximum.   Certainly, the oscillation time generalizes the period of solutions for the undamped model ($\alpha = 0$). 
In this investigation we analyze  the dependence of $\tau$ on $x_0$ and on $\alpha$ under the following working hypothesis:  
\begin{assumption}
\label{main_assumption}
On small $\epsilon-$neighborhood of $0$ the function $f$ is even and  for some constant $a>0$ we have
\[  
f(x) =  -a\, x^2 + O\left(\vert x\vert^4\right),  
\]
\end{assumption}

\begin{figure}
\centering
\includegraphics[width=0.75\textwidth]{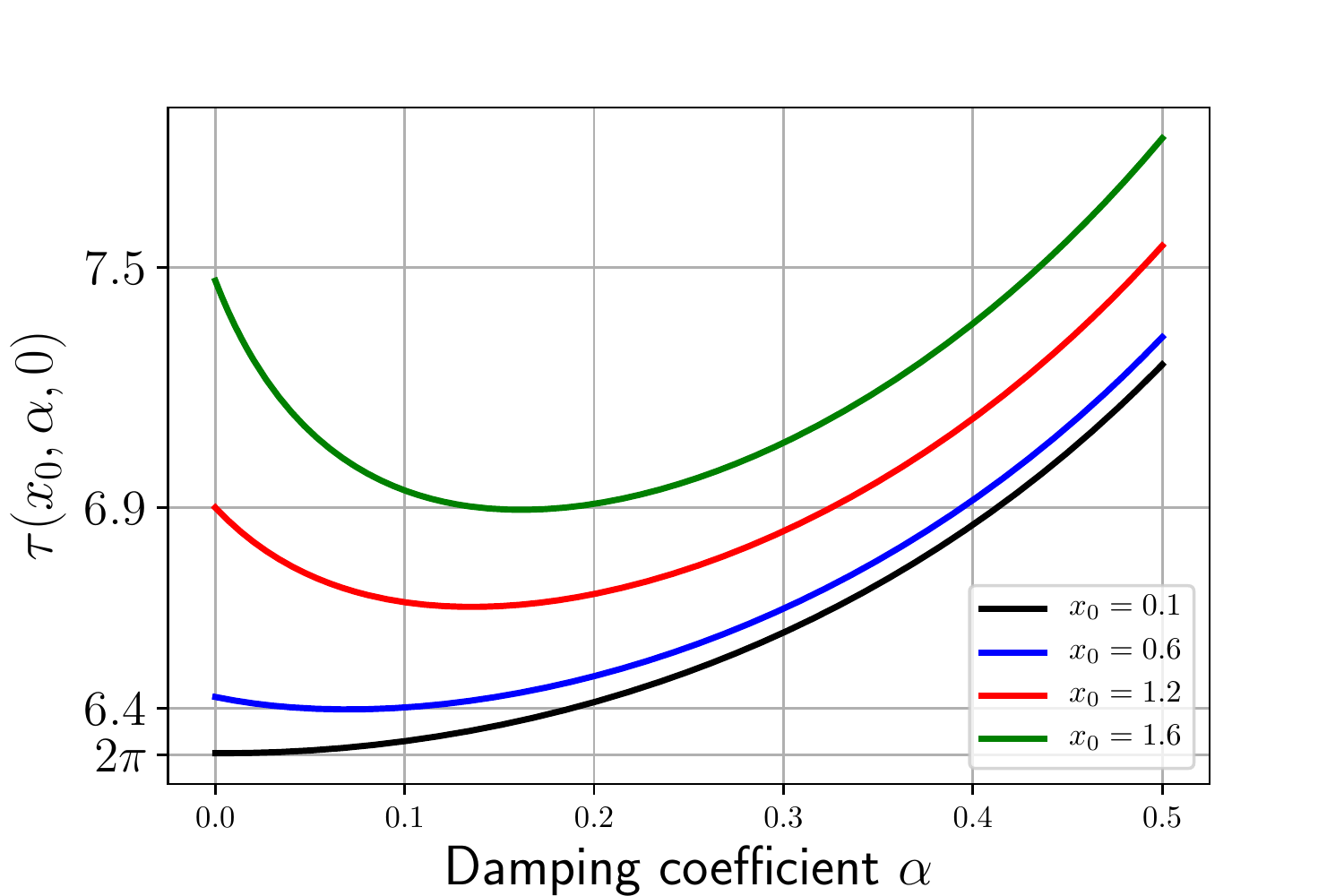}
\caption{\label{fig:tau}Numerical simulation  $\tau(x_0,\alpha)$ depending on $\alpha$ for several values of $x_0$. The nonlinear term $f$ was chosen so that $x(1+f(x)) = \sin{x}.$}
\end{figure}

We shall show  that for   $x_0$  fixed,   $\tau$ reaches a positive minimum at some $0<\alpha_0<1.$  It does not seem obvious that an increase in the damping coefficient $\alpha$ might cause a decrease in $\tau$. It is also worth noticing that the existence of a minimum of $\tau$ is a consequence  the sign of the constant $a$ in the above assumption.  Indeed, according to numerical experiments carried out by the author,  $\tau$ does not reach a positive minimum if $a<0$.  The author  is not aware of a similar result in the current literature nor whether this phenomena has been experimentally addressed.    The whole paper was written with the aim at the mathematical pendulum $x(1+f(x)) = \sin{x}.$  In that case, Figure \ref{fig:tau} summarize our findings by  picturing the numerically simulated value for $\tau(x_0,\alpha)$. Interestingly,  our qualitative analysis accurately reflects   variations of $\tau$  that are not easy to spot numerically. For instance, the minimum of $\tau(x_0, \alpha)$ for $x_0 = 0.1$ is not evident in Figure \ref{fig:tau}.

The arguments and proofs in this paper are entirely based on well established  techniques of ODE theory. However, the main result (Theorem \ref{thm:main}) rests on   delicate estimates  involving a differential equation describing the dependence of the solution $x=x(t)$ with respect to $\alpha$.

%%%%%%%%%%%%%%%%%%%%%%%%%%%%%%%%%%%%%%%%%%%%%%%%%%%%%%%%%%
\section{Underdamped oscillations}
Definitions of underdamped oscillations in linear systems naturally carry over to solutions of \eqref{eq:pendulum}. 
From now on, $x(\cdot, x_0, \alpha)$ stands for the unique solution to \eqref{eq:pendulum} satisfying the initial condition $x(0) = x_0$ and $\dot{x}(0) = 0.$ We also write $\tau(x_0,\alpha)$ to highlight the dependence of the oscillation time on $x_0$ and $\alpha$.  We will write simply $\tau$ or $x$ when no confusion can arise.    It is convenient to represent \eqref{eq:pendulum} in the phase space $(x,v)$ with $\dot{x}= v:$   
\begin{equation}
\label{eq:pendulum-phase}
\begin{split}
\dot{x}& =v \\
\dot{v} & =  -2\alpha\, v  -x  -x\,f\left(x\right).\\
\end{split}
\end{equation}
Equation \eqref{eq:pendulum-phase} is explicitly solvable whenever $f\equiv 0$, and in that case,  its solution is given by 
\begin{equation}
\label{eq:linear_xv}
\begin{split}
x_{l}(t) = & \frac{e^{-\alpha\,t}}{\omega}
\left( \omega\, \cos{\omega\, t} + \alpha\, \sin{\omega\, t} \right)x_0 \\
v_{l}(t) = & -\frac{e^{-\alpha\,t}}{\omega}\, \sin{\omega\, t}\, x_0
\end{split}
\end{equation}
where $\omega = \sqrt{1-\alpha^2}$. Moreover, the oscillation time $\tau_l$ is given by
\[
\tau_l = \frac{2\pi}{\omega} = \frac{2\pi}{\sqrt{1-\alpha^2}}.
\]
Notice that $\tau_l$ is an increasing function that solely depends on  $\alpha$.

Though a closed-form solution of  \eqref{eq:pendulum} is either not known or impractical, we could express the relevant solutions implicitly.  To that end,  we rewrite \eqref{eq:pendulum-phase} so that the nonlinear term $-x\,f\left(x\right)$ assumes the role of a non homogeneous forcing term. 
The expression for the solution $(x,v)$ is implicitly given by
\begin{equation}
\label{eq:var-par1}
\begin{split}
x(t) = &  x_{l}(t) - \frac{1}{\omega}  \int_{0}^{t}
e^{-\alpha\left(t-s\right)}\sin\omega\left( t-s\right)x(s)\, f(x(s))\, ds \\
v(t) = & v_{l}(t) - \frac{1}{\omega}  \int_{0}^{t}
e^{-\alpha\left(t-s\right)}\left(\omega\, \cos{\omega\left( t-s\right)}  -\alpha \sin\omega\left( t-s\right)\right)x(s)\, f(x(s))\, ds 
\end{split} 
\end{equation}

 Next, we estimate the solutions of \eqref{eq:pendulum-phase} in the conservative case ($\alpha=0$) in which all solutions are periodic and the period is given by $\tau \equiv \tau(x_0,0).$
\begin{lemma}
\label{lemm:estimate}
 If $(x,v)$ stands for the solution to \eqref{eq:pendulum-phase} with $\alpha = 0$  that satisfies $(x(0),v(0))= (x_0,0)$, then there exists $\delta>0$ so that for all $|x_0|\le \delta$ and all $0\le t \le \tau$ we have
\begin{equation}
 \label{eq:xv_conserv}
 x(t) = x_0\, \cos t + R_1(t, x_0), \quad 
 v(t) = -x_0\, \sin t + R_2(t, x_0), 
\end{equation}
where 
\[
|R_i(t, x_0)| \le \text{const }\, |x_0^3|, \quad i = 1,2.
\]
\begin{proof}
Letting $\alpha= 0$ in  \eqref{eq:var-par1} we obtain
\begin{equation}
\label{eq:help}
   R_1(t, x_0)=   -\int_{0}^{t}\cos(t-s)\, x(s)\, f(x(s))\, ds. 
\end{equation}
Since $(0,0)$ is a stable equilibrium solution to  \eqref{eq:pendulum-phase}, there exists $\delta>0$ and $\epsilon >0$ so that any solution $(x,v)$ to  \eqref{eq:pendulum-phase} starting at $(x_0,0)$, with $|x_0| \le \delta $ satisfies  $|x(t)| \le \epsilon$. 
Now write $F(z)= -z\, f(z)$ and notice that for some $\xi \in (-\epsilon, \epsilon)$ we have
\[
F(x(s)) = F\left(x_0 \cos{s}   + R_1\left(s, x_0\right)\right) 
 =  F\left(x_0 \cos{s}\right) + R_1\left(s, x_0\right)F'\left(\xi\right).
 \]
 Next, identity \eqref{eq:help}, Assumption \eqref{main_assumption} and some standard estimations yield
\[
|R_1(t, x_0)|\le 2a|x_0^3|+ c_2 \int_{0}^{t}  |R_1(s,x_0)|\,ds
\]
where $c_2 = \max_{z\in [-\epsilon, \epsilon]} |F'(z)|$. The first claim follows now from Gronwall's inequality. The proof of  the estimation for $R_2$ is analogous.
\end{proof}
\end{lemma}

At this point it is appropriated to define the \emph{half oscillation time} $\hat{\tau}= \hat{\tau}(x_0, \alpha)$ to be the time spent by the solution $x(t, x_0, \alpha),$ $t\ge 0,$ reaching the next local minimum. If $\alpha = 0$ and $f$ is even, the symmetry of the solution \eqref{eq:pendulum} yields. $2\hat{\tau} = \tau.$

\begin{lemma}
\label{lemm: half-tau}
If  $\hat{\tau}= \hat{\tau}(x_0, \alpha)$ denote the half oscillation time and $a$ is the constant of Assumption \ref{main_assumption}, then 
\[
\hat{\tau}(x_0, \alpha) >\frac{\pi}{\sqrt{1 - \alpha^2} } \quad \text{and} \quad \lim_{x_0 \to 0^{+}}\hat{\tau}(x_0, 0) =  \pi + \frac{a\, \pi}{8}\, x_0^2  + o(x_0^3). 
\]
\begin{proof}
We  introduce  introduce the polar coordinates
\[
r  =\sqrt{x^2+ v^2}, \quad  
\tan{\theta}  =  \frac{x}{v},
\]
to obtain 
\begin{equation}
\label{eq:r-theta}
\begin{split}
\dot{\theta}& =- \left( 1+  \alpha\,\sin{2\theta} + \sin^2{\theta} \, f\left(x\right)\right) \\
\dot{r} & = -\frac{v}{r}
 \left( 2\alpha \,   v + x\, f\left(x\right)\right)\\
\end{split}
\end{equation}

As a consequence of equation \eqref{eq:r-theta} we obtain the following expression for the half oscillation time $\hat{\tau}= \hat{\tau}(x_0, \alpha)$

\begin{equation}
\label{eq:tau}
\hat{\tau} = \int_{0}^{\pi}\frac{d\theta}{1+  \alpha\,\sin{2\theta} + \sin^2{\theta}\, f\left(x(\theta)\right) }.
\end{equation}
Now, the effect of the nonlinearity on the oscillation time is clear. By Assumption \ref{main_assumption} we obtain
\[  
\hat{\tau}(x_0, \alpha) > \int_{0}^{\pi}\frac{d\theta }{1+  \alpha\,\sin{2\theta}  }= \frac{\pi}{\sqrt{1 - \alpha^2} }.
\]

For $\alpha = 0$ we use estimation  \eqref{eq:xv_conserv} to obtain
\[  
\hat{\tau}(x_0, 0) = \int_{0}^{\pi}\frac{d\theta}{1    - a\,x_0^2\, \sin^2{\theta}\, \cos^2{t(\theta)} }  + o(x_0^3).
\]
Now a straightforwards computation yields
\[ 
\lim_{x_0 \to 0^{+}}\hat{\tau}(x_0, 0) = \pi,  \; 
\lim_{x_0 \to 0^{+}}\frac{\partial \hat{\tau}}{\partial x_0}(x_0, 0) = 0.
\]
Now, the expression for $\frac{\partial^2 \hat{\tau}}{\partial x_0^2}(x_0, 0)$ is somewhat cumbersome. However,  taking into account that $ \lim_{x_0 \to 0} t(\theta)= \theta$, we readily obtain 
\[
\lim_{x_0 \to 0^{+}}\frac{\partial^2 \hat{\tau}}{\partial x_0^2}(x_0, 0) = \int_{0}^{\pi} 2a\, \sin^2{\theta}\, \cos^2{\theta}\, d\theta=\frac{2\, a\pi}{8}, 
\]
and the second claim of the lemma follows by the second order Taylor expansion of $\hat{\tau}(x_0, 0)$ around $x_0$
\end{proof}
\end{lemma}

A reasoning analogous to that in the proof of the preceding lemma shows that 
\[
 \tau(x_0, \alpha) >\frac{2\pi}{\sqrt{1 - \alpha^2} }\equiv \tau_l.
\] This inequality is illustrated in Figure  \ref{fig:tau_comp} when $a = 1$. Had we considered in Assumption \ref{main_assumption}   negative values for $a$, then the inequality would reverse to $\tau(x_0, \alpha) < \tau_l$ as it is depicted in Figure  \ref{fig:tau_comp}.
%%%%%%%%%%%%%%%%%%%%%%%%%%%%%%%%%%%%%%%%%%
%%%%%%%%%%%%%%%%%%%%%%%%%%%%%%%%%%%%%%%%%%%%%%%%%%%%%
\section{The role of the viscous damping}
It is not difficult at all to obtain a differential equation describing the  movement of the pendulum depending on the viscous damping coefficient. Indeed, writing 
\[
X(t, x_0, \alpha) = \frac{\partial x}{\partial \alpha}\left(t,x_0,  \alpha \right), \quad 
V(t, x_0, \alpha) = \frac{\partial v}{\partial \alpha}\left(t,x_0,  \alpha \right).
\]
Derivation of equation  \eqref{eq:pendulum-phase} with respect to $\alpha$ yields:
\begin{equation}
\label{eq:XV}
\begin{split}
\dot{X} & = V\\
\dot{V} & =- 2\alpha\, V -X - 2v -
 \left(  x\,  f'\left(x\right)+ f\left(x\right) 
\right)X.
\end{split}
\end{equation}
As for the initial conditions we have 
 \[
X(0, x_0, \alpha) = 0, \quad 
V(0, x_0, \alpha) = 0.
\]
Let us write $G(x) = -\frac{d}{dx }\left(x\, f\left(x\right)\right)$.  
Again, as we did with equation \eqref{eq:pendulum-phase}, equation \eqref{eq:XV}   can be seen as a linear homogeneous part plus the forcing term $- 2v +G(x)\,X.$ 
The solution $X,V$ is implicitly given by
\[  
\begin{split}
X(t) = &     \frac{1}{\omega}  \int_{0}^{t} 
e^{-\alpha\left(t-s\right)}\sin\omega\left( t-s\right)\big\{ \\
  & \qquad - 2v(s)+G(x(s))\, X(s)\left. \right\} ds \\
V(t) = &   \frac{1}{\omega}  \int_{0}^{t}
e^{-\alpha\left(t-s\right)}\left(\omega\, \cos{\omega\left( t-s\right)}  -\alpha \sin\omega\left( t-s\right)\right)\big\{\\
 & \qquad  - 2v(s)+G(x(s))\, X(s)\left. \right\} ds
\end{split} 
\]
In particular, for $\alpha = 0$ the above expressions reduce to 
\begin{equation}
\label{eqn:var-par2}
\begin{split}
X(t) = &        \int_{0}^{t} 
\sin \left( t-s\right)\big\{- 2v(s)+G(x(s))\, X(s) \big\}\, ds \\
V(t) = &\int_{0}^{t} \cos{ \left( t-s\right)} \big\{- 2v(s)+G(x(s))\, X(s) \big\} \, ds
\end{split} 
\end{equation}

The following lemma does  the heavy lifting to deliver the main result of the paper. 

\begin{lemma}
\label{lemm: heavy-lifting}
Under Assumption \ref{main_assumption}, if $\hat{\tau}=\hat{\tau}(x_0, 0) $ denotes the half oscillation time when $\alpha = 0$, then for $0<x_0\ll 1$ we have
$
V(\hat{\tau}, x_0, 0) >0.
$
\begin{proof} We start with an auxiliary estimate for $X(t)$ in equation \eqref{eqn:var-par2}. By  Lemma \eqref{lemm:estimate} and by Assumption \ref{main_assumption}, for $0<t \le \pi$ we have
\begin{equation}
\label{eq:X}
 X(t) = x_0\left(-t\, \cos{t}+ \sin{t}\right) +    3a\, x_0^2 \,\int_{0}^{t} 
\sin \left( t-s\right)\cos^2{s} \, X(s) \, ds + O(|x_0|^4)   
\end{equation}
Notice that $X_1(t) \equiv x_0\left(-t\, \cos{t}+ \sin{t}\right)$ does not vanish on $(0, \pi)$ and that $G(x(s))>0$ provided $0<x_0\ll 1$.  Further, the initial conditions for $X(t)$ at $t= 0$ and equation \eqref{eq:XV} yield that
\[
X(0) = 0 = \dot{X}(0) = \ddot{X}(0)\text{ and } \dddot{X}(0)= 2x_0\left( 1 + f\left(x_0\right) \right)>0,
\]
meaning that $X(t)$ is positive on an interval  $(0, \epsilon)$ with $\epsilon>0$. We claim that $X(t)>0$ for $0<t \le \pi.$ On the contrary, there exists  $\epsilon <t_0<\pi$  such that $X(t_0) = 0$ and $X(t)>0 $ for $t\in (0, t_0)$. Now, by Lemma \ref{lemm: half-tau} we know that 
$\hat{\tau}> \pi$. Therefore, the polar angle $\theta(t) $ in  \eqref{eq:r-theta} satisfies $-\pi< \theta(t) <0$ for all $0<t<\pi$ and a fortiori  $v(t)<0$ on $(0, \pi]$. But this is a contradiction to the first equation of \eqref{eqn:var-par2} evaluated at $t= t_0$ since for $s\in (0, t_0)$ we have
\[
\sin \left( t_0-s\right)\big\{- 2v(s)+G(x(s))\, X(s) \big\} >0.
\]

Next, by the  equation \eqref{eq:X} it follows immediately that 
$X(t) =X_1(t) +  O(|x_0|^3)$. Analogously, for $V(t)$ we obtain 
\[
\begin{split}
V(t) = & x_0\, t\, \sin{t} + 3a\, x_0^2   \int_{0}^{t} 
\cos \left( t-s\right)\cos^2{s} \, X_1(s) \, ds + O(|x_0|^4)\\
\equiv  & V_1(t) + V_2(t)+ O(|x_0|^4)
\end{split}
\]
where $V_1(t) \equiv x_0\, t\, \sin{t}$. Now, $V_2(t)$ can be explicitly evaluated. For the reader's convenience,  we write the complete expression for $V_2:$ 
\[
\begin{split}
 V_2(t) = & 3a\, x_0^3\Big(-\frac{1}{32}\,(6\,t^2 + 5)\,\cos{t} - \frac{3}{32}\,t\,\sin{3\,t} -  \\
  & \frac{1}{16}\,t\,\sin{t} - \frac{17}{128}\,\cos{3\,t} + \frac{37}{128}\,\cos{t}\Big).
\end{split}
\]
Moreover, it is somewhat tedious but straightforward to show that $V_2$ is positive and  increasing on a small neighborhood of $\pi$.  By Lemma \ref{lemm: half-tau} $\hat{\tau}>\pi$, therefore  
\[ 
V_2(\hat{\tau})> V_2(\pi)= \frac{9a\, x_0^3\, \pi^2}{16}.
\]
Again, by Lemma \ref{lemm: half-tau} we obtain
\[
\begin{split}
V_1(\hat{\tau}) = & V_1(\pi) + \left(\hat{\tau}-\pi\right)V_1' (\pi)+  O(|x_0|^4)\\
= & - \frac{a\, x_0^3\, \pi^2}{8}+  O(|x_0|^4), 
\end{split}
\]
so that $V(\hat{\tau}) = V_1(\hat{\tau})+ V_2(\hat{\tau})>0.$
\end{proof}
\end{lemma}

Now we are in a position to show the main result of the paper
\begin{theorem}
\label{thm:main}
Under Assumption \ref{main_assumption}, there exists a $\delta >0$ such that for $0<x_0<\delta$ fixed, the oscillation time $\tau(x_0, \alpha)$, for $0 < \alpha<1$,  reaches a positive minimum at some   $0 < \alpha<1$. Moreover, 
\[
\lim_{ \alpha \to 1^{-}}\tau(x_0, \alpha)= \infty.
\]
\begin{proof}
We let $0<x_0\ll 1$ fixed by now and denote by  $(x,v)$ be the solution of equation  \eqref{eq:pendulum-phase}. By definition of $\hat{\tau}$ we have $v(\hat{\tau}, \alpha) = 0$, so that   the Implicit Function Theorem yields 
\[
\frac{\partial \hat{\tau}}{\partial \alpha}\,\dot{v}(\hat{\tau},   \alpha) + V(\hat{\tau},   \alpha)=0,
\] therefore
\[
\frac{\partial \hat{\tau}}{\partial \alpha} = 
\frac{V(\hat{\tau}, \alpha)}
{x\left(\hat{\tau}, \alpha \right)\left( 1+f\left(  x\left(\hat{\tau}, \alpha \right)\right)\right) }.
\]
Since 
$x\left(\hat{\tau}, \alpha \right)$ is negative, it follows from  Lemma \ref{lemm: heavy-lifting} that and $\frac{\partial \hat{\tau}}{\partial \alpha}\vert_{\alpha= 0}<0$. Now we shall show that the last inequality holds for the oscillation time $\tau$. To do that, write $\hat{x}_0 = -x(\hat{\tau}(\alpha, x_0), x_0)$ and see that 
\[
\tau(\alpha, x_0) = \hat{\tau}(\alpha, x_0)+  \hat{\tau}(\alpha, \hat{x}_0).
\]
That is to say, the half oscillation time depends on $|x_0|$ only. Notice that $\hat{x}_0\le x_0$  and the equality holds in the conservative case $\alpha= 0$ only. Therefore
\[
\frac{\partial  \tau}{\partial \alpha}(\alpha, x_0) = \frac{\partial \hat{\tau}}{\partial \alpha}(\alpha, x_0)  +\frac{\partial \hat{\tau}}{\partial \alpha}(\alpha, \hat{x}_0)  - \frac{\partial  \hat{x}_0}{\partial \alpha}(\alpha, x_0)\, \frac{\partial \hat{\tau}}{\partial \alpha}(\alpha, x_0) = 0. 
\]
Moreover, since 
\[
\frac{\partial  \hat{x}_0}{\partial \alpha}(\alpha, x_0)= v(\hat{\tau}(\alpha, x_0), x_0)=0,
\]
we have that 
\[
\lim_{x_0 \to 0^{+}} \frac{\partial  \tau}{\partial \alpha}(\alpha, x_0) = 2\lim_{x_0 \to 0^{+}} \frac{\partial  \hat{x}_0}{\partial \alpha}(\alpha, x_0)
\]
Finally, by the first claim of Lemma \ref{lemm: half-tau}, $ \tau(\alpha, x_0)$ must attain a minimum at some $0<\alpha<1$. 
\end{proof}
\end{theorem}
%%%%%%%%%%%%%%%%%%%%%%%%%%%%%%%%%%%%%%%%%%%%%%%%%%%%%%%
%%%%%%%%%%%%%%%%%%%%%%%%%%%%%%%%%%%%%
\section{Conclusions and final remarks} 

An oscillating mass exhibits gradually diminishing amplitude in the presence of damping. The time spent by the mass completing one oscillation depends on several factors, as the model for the restoring force, how the oscillation starts, and the nature of the damping. For the sake of our discussion we consider a vertical pendulum with a nonlinear  restoring force resembling the mathematical pendulum, letting the oscillation start at a small amplitude with vanishing velocity and a viscous damping model with a (normalized) viscosity coefficient $\alpha$.  We have proved that the oscillation time $\tau \equiv \tau(\alpha)$ does not depend monotonically on $\alpha$, meaning that there exists a threshold $\alpha_0$ (which depends on the starting amplitude of the oscillation) such that $\tau$ reaches a local minimum  at $\alpha_0$ (see Figure \ref{fig:tau}). It is worth noticing that this behavior cannot be observed if the restitution force is linear, i.\ e., what we report in this paper is essentially a nonlinear phenomenon.    

\begin{figure}
\centering
\includegraphics[width=0.75\textwidth]{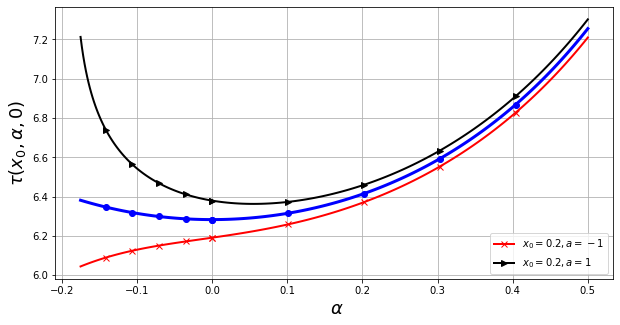}
\caption{\label{fig:tau_comp} Numerical simulation  of the oscillation time $\tau$ depending on the damping coefficient $\alpha$ with starting amplitude  $x_0 = 0.2$ and non linear restoring term given by  $f(x)= -a\, x^2$, $a =\pm 1$. The curve with the round marker (blue in the online version) corresponds to the oscillation time $\tau_l$ of the linear case $f\equiv 0$}
\end{figure}

The proof of existence of a positive minimum for the oscillation time rests heavily on the fact that the constant $a$ in Assumption \ref{main_assumption} is positive. Just to experiment the effect of changing the sign of the constant $a$,  we carried out some numerical simulations of $\tau$ with the nonlinear term $f(x) = -a\, x^2$ for $a = 1, -1$. The corresponding equations are  particular cases of an \textit{unforced Duffing oscillator} \cite{wiggins}. The numerical results are shown in Figure \ref{fig:tau_comp}. Just for the sake of the numerical experimentation we also  considered negative values for $\alpha$. If $a = 1$ we see that $\tau$ reaches its minimum at a positive value for $\alpha$.  By contrast, if $a = -1$ no minimum seems to exist.      The curve with the round marker (blue in the online version) corresponds to the oscillation time of the linear case $\tau_l = \frac{2\pi}{\sqrt{1-\alpha^2}}$. 

The numerical experimentation of the oscillation time $\tau$ (not shown in this paper)   assuming a quadratic damping exhibits the same behavior as the graphics of Figure \ref{fig:tau_comp}. If the readers are  curious about the numerical experiments,   they could take a look at the author's GitHub page
\begin{center}
\url{https://github.com/arangogithub/Oscillation-time}
\end{center}
and download a Jupyter notebook with the python code featuring the results shown in Figures \ref{fig:tau} and \ref{fig:tau_comp}

 \section*{Acknowledgment} The author  would like to give the reviewer his very heartfelt thanks for carefully reading the
manuscript and for pointing out several inaccuracies of the document.

\bibliographystyle{plain}
\bibliography{main}

\end{document}